\documentclass{article}
\usepackage{graphicx}
\usepackage{amsmath}

\begin{document}

\title{The representations of polyadic-like equality algebras}
\author{Mikl\'{o}s Ferenczi}
\maketitle

\begin{abstract}
It is proven that Boolean set algebras with unit $V$ of the form $%
\bigcup\limits_{k\in K}$ $^{\alpha }U_{k}^{{}}$ are axiomatizable (i.e. if $%
V $ is a union of Cartesian products). The axiomatization coincides with
that of \ cylindric polyadic equality algebras (class CPE$_{\alpha }$). This
is an algebraic representation theorem for the class CPE$_{\alpha }$ by
relativized polyadic set algebras in the class Gp$_{\alpha } $$. $ Similar
representation theorems are claimed for the classes strong cylindric
polyadic equality algebras (CPES$_{\alpha }$) and cylindric $m$-quasi
polyadic equality algebras ($_{m}$CPE$_{\alpha }$). These are polyadic-like
equality algebras with infinite substitution operators and single
cylindrifications. They can be regarded also as infinite transformation
systems equipped with diagonals and cylindrifications. No representation
theorem or neat embedding theorem has proven for this class of algebras yet,
except for the locally finite case. The theorems proven here answer some
unsolved problems.

subjclass: 03G15, 03G25

keywords: polyadic algebras, cylindric algebras
\end{abstract}

\bigskip \bigskip

This is a shortened version of the paper submitted, with the same title. In
\cite{FePol} it is shown that the class of Boolean set algebras with unit $%
V, $\ such that $V$\ is a union of weak Cartesian spaces, i.e.,\textit{\ }$%
V= $\textit{\ }$\bigcup\limits_{k\in K}$\textit{\ }$^{\alpha }U_{k}^{(pk)}$%
\textit{\ }is axiomatizable (class Gwp$_{\alpha }$). The axiomatization thus
obtained coincides with that of\textit{\ }the transposition algebras (TA$%
_{\alpha }$). Transposition algebras are Boolean algebras extended by the
abstract transposition operators ($p_{ij})$, single substitutions ($%
s_{j}^{i} $), cylindrifications ($c_{i}$) and the diagonal constants ($%
d_{ij} $), where $i,j\in \alpha $. These algebras are definitionally
equivalent to cylindric quasi polyadic equality algebras (CQE$_{\alpha }$),
i.e., to Boolean algebras with \textit{finite} substitutions ($s_{\tau }$,
where $\tau $ is a\textit{\ finite} transformation on $\alpha $),
cylindrifications ($c_{i}$) and diagonal constants ($d_{ij}$), where $i,j\in
\alpha $.\newline
In this paper, it is shown, among others, that the class of Boolean set
algebras with unit $V,$\ such that $V$\ is a union of Cartesian spaces,
i.e., $V=$\ $\bigcup\limits_{k\in K}$\ $^{\alpha }U_{k}^{{}}$ is \textit{%
axiomatizable} (class Gp$_{\alpha }$). The resulting axiomatization
coincides with that of the cylindric polyadic equality algebras (class CPE$%
_{\alpha }$), see the main representation theorem, Theorem 3. These algebras
are Boolean algebras extended by substitutions ($s_{\tau },$ where $\tau \in
$ $^{\alpha }\alpha $ is \textit{arbitrary}), cylindrifications ($c_{i}$)
and the diagonal constants ($d_{ij}$), where $i,j\in \alpha $. The attribute
``cylindric'' which qualifies the class CPE$_{\alpha }$ comes from the fact
that only single cylindrifications are allowed here (in contrast with
polyadic algebras where the operation ``quantifier $\exists $'' is defined
for any subset of $\alpha $). This is why, among others, the algebras
occuring here are considered as \textit{polyadic-like} algebras. The
difference between the classes CQE$_{\alpha }$ and CPE$_{\alpha }$ is in the
\emph{finite} or \emph{infinite} nature of the transformations occuring in
the substitutions.\newline
The techniques used in \cite{FePol} and here are essentially \emph{different}%
. While in \cite{FePol} finite methods are used (which cannot be used here),
now, the neat embedding technique is applied (which cannot be applied to the
algebras in CQE$_{\alpha }$).\newline
Up until now, we have looked at our subject from the viewpoint of the
axiomatizability of certain classes of set algebras (Gwp$_{\alpha }$ and Gp$%
_{\alpha }$), i.e., from the viewpoint of set theory and\ geometry\emph{.}
Now, let us consider our subject from the viewpoint of the representation
theory of\ algebraic logic.\newline
The quasi polyadic (CQE$_{\alpha }$) and the cylindric polyadic (CPE$%
_{\alpha }$) cases (where the substitutions are finite and infinite,
respectively) share several common features besides differing with respect
to others.\newline
As for the common properties, neither CQE$_{\alpha }$ nor CPE$_{\alpha }$
can be represented in the classical sense (i.e., as subdirect product of
polyadic set algebras). Also, in both cases, the representability given here
is achieved by\textit{\ relativized} set algebras (Gwp$_{\alpha }$ and Gp$%
_{\alpha }$, respectively), instead of ordinary set algebras. It is also
important to emphasize that the commutativity of cylindrifications is not
assumed in these classes, either.\newline
Some words about the \emph{neat embedding technique} applied in the paper.
The proof of the main representation theorem (Theorem 3) follows a classical
line of thought. The proof consist of two steps. The first step is based on
a theorem which establishes that the algebra in question is neatly
embeddable (now, the celebrated theorem of Daigneault$,$ Monk and Keisler
and its versions are used, see \cite{He-Mo-Ta}, II. Thm. 5.4.17). The second
step of the proof is based on the \textit{neat embedding theorem} Theorem 1
proven here. A neat embedding theorem for a given class says that neatly
embeddability is equivalent to representability. With this kinds of proofs
of the representation theorems, the \emph{Henkin style proofs} of
completeness theorems can be associated in mathematical logic.\newline
A few words about the \textit{history} of the topic may be in order here.
Daigneault, Monk and Keisler proved that polyadic algebras (without
equality) are representable in the classical sense, but, polyadic \textit{%
equality} algebras, in general, are not (\cite{Da-Mo}, \cite{Ke}). For
polyadic equality algebras and cylindric algebras, the representability was
proved only for special, but important classes, e.g., for the class of
locally finite algebras (see \cite{Ha56}, \cite{Ha57}, \cite{He-Mo-Ta}).%
\newline
\textit{Relativized set algebras }were used for representation first in the
theory of cylindric algebras. Using an idea of Henkin, Resek and Thompson
showed that extending the cylindric axioms by the so called merry-go-round
property, cylindric algebras become to be representable by cylindric
relativized set algebras. First, Andr\'{e}ka and Thompson published a proof
for this theorem (see \cite{An-Th}), later, the present author gave a
simpler axiomatization for the class in question (see \cite{Fe07} and \cite
{Fe-Sa}). The method used in Andr\'{e}ka and Thompson's proof is detailed in
\cite{Hi-HoStep}.\newline
The $r$-representation of \textit{polyadic} equality algebras by relativized
set algebras was investigated in \cite{FeMLQ}, \cite{FePol} and \cite{A-F-N}%
. In \cite{FeMLQ} it is shown that the background of the merry-go-round
property is an axiom of transposition algebras (and it is also a property of
polyadic algebras) and Resek and Thompson's theorem is closely related to
transposition algebras. The direct predecessor of the research here is the
investigation of the representability of transposition algebras in \cite
{FePol}.

We only deal with infinite dimensional algebras because the finite case is
an instance of the quasi polyadic case (see \cite{FePol}). Theorem 1 is a
neat embedding theorem, it is the key result to Theorem 3. Theorem 3 states
the $r$-representability of the classes cylindric polyadic equality algebras
(for CPE$_{\alpha }$) and strong cylindric polyadic equality algebras (class
CPES$_{\alpha }$). This theorem answers affirmatively two classical problems
raised in the literature. The one is the representability of infinite
transformation systems equipped with diagonals and cylindrifications (see
\cite{Slom} and \cite{Da-Mo}), and the other is whether the class Gp$%
_{\alpha }$ is a \textit{variety.} Theorem 4 is a consequence the
representation theorems.

\bigskip

\bigskip Now we list the important results with the necessary concepts:

\bigskip The concept of polyadic equality algebra is assumed to be known
(see \cite{He-Mo-Ta}, 5.4.1)

\bigskip

\textbf{Definition} \textbf{1} \textit{Assume that }$m<\alpha $\textit{\ and
}$m$\textit{\ is infinite}$.$\textit{\ Given a set }$U$\textit{\ and a fixed
sequence }$p\in $\textit{\ }$^{\alpha }U$\textit{, the set }
\begin{equation*}
_{m}^{\alpha }U_{{}}^{(p)}=\left\{ x\in \;^{\alpha }U\;:\;x\;\text{and\ }p\;%
\text{are different at most in }\mathit{m}\text{ many places}\right\}
\end{equation*}
\textit{is called the }$m$\textit{-}weak space\textit{\emph{\ }(or }$m$%
\textit{-weak Cartesian space) determined by }$p$\textit{\ and the base }$U$%
\textit{.}

\bigskip

\bigskip

Recall that the finite version ($m$ is finite) of the above definition is
the concept of the weak space, in notation $^{\alpha }U_{{}}^{(p)}$ (\cite
{He-Mo-Ta}, II. 3.1.2). A unit $V$ of some $\mathcal{A}\in $ Cprs$_{\alpha }$
can be composed as a disjoint union of certain subsets of $m$-weak spaces.
These latter subsets are called $m$-\textit{subunits} of $\mathcal{A}$, the
bases of these $m$-subunits are called $m$-\emph{subbases}\textit{\ }of $%
\mathcal{A}$\textit{.}

\bigskip

\bigskip

\textbf{Definition 2} (class $_{m}$Gwp$_{\alpha }$) \textit{A set algebra in
}Cprs$_{\alpha }$\textit{\ is called a} generalized $m$-weak polyadic
relativized set algebra\ \textit{if there are sets }$U_{k}$\textit{\ and }$%
p_{k}\in $ $^{\alpha }U_{k}$ \textit{such that }$V=\bigcup\nolimits_{k\in
K}{}$ $_{m}^{\alpha }U_{k}^{(p_{k})},$\textit{\ where }$V$\textit{\ is the
unit. The subclass of }$_{m}$Gwp$_{\alpha }$ \textit{such that the
disjointness of the }$U_{k}$'\textit{s is assumed is denoted by }$\overset{%
\bullet }{_{m}\text{Gwp}_{\alpha }}.$

\bigskip

\textbf{Definition 3} (class Gp$_{\alpha }$) \textit{A set algebra in }Cprs$%
_{\alpha }$\textit{\ is called a} generalized polyadic relativized set
algebra\ \textit{if there are sets }$U_{k}$\textit{\ such that }$%
V=\bigcup\nolimits_{k\in K}{}^{\alpha }U_{k},$\textit{\ where }$V$\textit{\
is the unit. The subclass of }Gp$_{\alpha }$ \textit{such that the
disjointness of the }$U_{k}$'\textit{s is assumed is denoted by} $\overset{%
\bullet }{\text{Gp}_{\alpha }}.$

\bigskip

\textbf{Definition 4 }(CPE$_{\alpha }$) \textit{A} cylindric polyadic
equality algebra\emph{\ }\textit{is a polyadic algebra such that instead of
the cylindrifications }$c_{\Gamma }$ \textit{only single cylindrifications }$%
c_{i}$ \textit{are defined, the cylindrifications are non-commutative and
the following weakening of the axiom} (P\textit{11) is assumed: }$\ $(P11)$%
^{\ast }:c_{i}s_{\sigma }x\leq s_{\sigma }c_{j}x,$ if $\left\{ j\right\} $= $%
\sigma ^{-1\ast }\left\{ i\right\} ,$ $\left\{ j\right\} \neq \emptyset ,$
and $c_{i}s_{\sigma }x=s_{\sigma }x$ else.

\bigskip

\textbf{Definition} \textbf{5} (CPES$_{\alpha }$) \textit{The class of }%
strong cylindric polyadic algebras \textit{is a subclass of }CPE$_{\alpha }$
\textit{such that the cylindrifications are commutative and} (P11) \textit{%
is assumed for the single cylindrifications.}

\bigskip

\bigskip \textit{A transformation }$\tau $\textit{\ defined on }$\alpha $%
\textit{\ is said to be }an $m$- transformation\textit{\ if }$\tau i=i$%
\textit{\ except for }$m$\textit{-many }$i\in \alpha .$\textit{\ The class
of }$m$-\textit{transformations is denoted by }$_{m}$T$_{\alpha }.$

\bigskip

\textbf{Definition 6} (class $_{m}$CPE$_{\alpha }$) \textit{A} cylindric $m$%
-quasi polyadic equality algebra \textit{of dimension }$\alpha $\textit{\ is
an algebra in }CPE$_{\alpha }$ \textit{such that the transformations }$\tau $
\textit{and} $\sigma $ \textit{are }$m$-\textit{transformation in the
definition of} CPE$_{\alpha }$\textit{, }i.e., $\tau ,\sigma \in $ $_{m}$T$%
_{\alpha }.$

\bigskip

\bigskip

The following one is a neat embedding theorem for the class $_{m}$CPE$%
_{\alpha }\cap $ L$m_{\alpha }$:

\bigskip

\textbf{Theorem 1 }\textit{If }$\mathcal{A}\in $ $_{m}$CPE$_{\alpha }\cap $ L%
$m_{\alpha },$ \textit{where }$m$\textit{\ is infinite, }$m<\alpha $\textit{,%
} then $\mathcal{A}\in $ SNr$_{\alpha }\mathcal{B}$\textit{\ for some} $%
\mathcal{B}$ $\in $ $_{m}$CPE$_{\alpha +\varepsilon },$ \textit{where }$%
\varepsilon $\textit{\ is infinite, if and only if\ }$\mathcal{A}\in $ Is$%
_{m}$Gwp$_{\alpha }.$

\bigskip

\bigskip The following theorem is the main representation theorem of the
paper:

\bigskip

\bigskip \textbf{Theorem 3} \ \textit{If} $\mathcal{A}\in $ CPE$_{\alpha
}\cup $ CPES$_{\alpha },$ \textit{then} $\mathcal{A}\in $ IsGp$_{\alpha }$.

\bigskip

\bigskip For the class CPE$_{\alpha },$ the converse of the theorem also
holds, i.e. $\mathcal{A}\in $ CPE$_{\alpha }$ if and only if $\mathcal{A}\in
$ IsGp$_{\alpha }.$

Now, we state again a \textit{neat embedding theorem}, a theorem for the
class CPES$_{\alpha }.$

\bigskip

\textbf{Theorem 4 }\textit{Assume that} $\mathcal{A}$ $\in $ CPES$_{\alpha }$
\textit{and }$\alpha $\textit{\ is infinite.Then, the following propositions
}(i) \textit{and }(ii) \textit{are equivalent:}

(i) $\mathcal{A}$ $\in $ SNr$_{\alpha }$CPES$_{\alpha +\varepsilon }$
\textit{for some infinite }$\varepsilon $

(ii) $\mathcal{A}\in $ Is(Gp$_{\alpha }$ $\cap $ Mod$\left\{ \text{(C4)},%
\text{ (P11)}^{\ast }\right\} $)

Ferenczi, M.

ferenczi@math.bme.hu


\bigskip

\begin{thebibliography}{99}
\bibitem{A-F-N}  Andr\'{e}ka, H., Ferenczi, M., N\'{e}meti, I., \textit{%
Cylindric-like Algebras and Algebraic Logic}, Springer, to appear

\bibitem{An-Th}  H.~Andr\'{e}ka and R.~J.~Thompson, \textit{A Stone type
representation theorem for algebras of relations of higher rank},
Transaction of Amer. Math. Soc., \textbf{309} (2) (1988), 671--682.

\bibitem{Da-Mo}  A.~Daigneault and J.~D.~ Monk, \textit{Representation
theory for polyadic algebras}, Fund. Math., \textbf{52} (1963), 151--176.

\bibitem{Fe-Sa}  M.~Ferenczi and G.~S\'{a}gi, \textit{On some developments
in the representation theory of cylindric-like algebras}, Algebra
Universalis, \textbf{55} (2--3) (2006), 345--353.

\bibitem{Fe07}  M.~Ferenczi, \textit{On cylindric algebras satisfying
merry-go-round properties}, Logic Journal of IGPL, \textbf{15} (2) (2007),
183--197.

\bibitem{FeMLQ}  M.~Ferenczi, \textit{Partial transposition implies
representability in cylindric algebras}, Mathematical Logic Quarterly,
\textbf{57}, 1, (2011), 87--94

\bibitem{FePol}  M. Ferenczi, \textit{The polyadic generalization of the
Boolean axiomatization of fields of sets, }Transaction of American Math.
Soc., to appear

\bibitem{Ha56}  P. Halmos, Algebraic logic II., \textit{Homogeneous, locally
finite polyadic Boolean algebras.} Fundamenta Mathematicae, \textbf{43},
(1956), 255--325

\bibitem{Ha57}  P. Halmos, Algebraic logic IV., \textit{Equality in polyadic
algebras}, Transaction of American Math. Soc.,\textbf{\ 86}, (1957), 1--27

\bibitem{He-Mo-Ta}  L.~Henkin, J.~D.~ Monk and A.~Tarski, \textit{Cylindric
Algebras I-II}, North Holland, 1985.

\bibitem{He-Mo-TaCikk}  L.~Henkin, J.~D.~ Monk and A.~Tarski, Representable
cylindric algebras, Annals of Pure and Applied Logic, 31, 10, (1986), 23-60

\bibitem{Hi-HoStep}  R. Hirsch and I. Hodkinson, \textit{Step by
step-building representations in algebraic logic,} J. Symbolic Logic Volume
62 (1997), 225-279.

\bibitem{Hi-Hokonyv}  R.~Hirsch and I.~Hodkinson, \textit{Relation Algebras
by Games}, North Holland, 2002.

\bibitem{Ke}  Keisler, H. J., \textit{A complete first order logic with
infiniteary predicates}, Fund. Math., \textbf{52} (1963), 177--203

\bibitem{Ta}  Tarek Sayed A., \textit{Some results about neat reducts},
Algebra Universalis, \textbf{1} (2010) 17--36

\bibitem{Takonyv}  Tarek Sayed A., \textit{Algebraic Logic}, manuscript
\end{thebibliography}
\end{document}